\font\script=eusm10.
\font\sets=msbm10.
\font\symbols=msam10.
\font\stampatello=cmcsc10.

\def\1{{\bf 1}}
\def\sgn{{\rm sgn}}
\def\oneH{{\1_{_H}}}
\def\wH{{w_{_H}}}
\def\FTwH{{\widehat{w}_{_H}}}
\def\wHone{{w'_{_H}}}
\def\wHtwo{{w''_{_H}}}
\def\uH{{u_{_H}}}
\def\FTuH{{\widehat{u}_{_H}}}
\def\sgnH{{\rm sgn}_{_H}}
\def\FTsgnH{{\widehat{\sgn}_{_H}}}
\def\CH{{C_{_H}}}
\def\FTCH{{\widehat{C}_{_H}}}
\def\avesum{\sum_{x\sim N}}
\def\dashsum{\mathop{\enspace{\sum}'}}
\def\starsum{\mathop{\enspace{\sum}^{\ast}}}
\def\square{\hbox{\vrule\vbox{\hrule\phantom{s}\hrule}\vrule}}
\def\defineq{\buildrel{def}\over{=}}
\def\defin{\buildrel{def}\over{\Longleftrightarrow}}
\def\doublesum{\mathop{\sum\sum}}
\def\multiplesum{\mathop{\sum \enspace \cdots \enspace \sum}}
\def\supporto{{\rm supp}\,}
\def\C{\hbox{\sets C}}
\def\N{\hbox{\sets N}}
\def\R{\hbox{\sets R}}
\def\Z{\hbox{\sets Z}}
\def\Corr{\hbox{\script C}}
\def\divisor{\hbox{\bf d}}
\def\Res{\mathop{{\rm Res}\,}}
\def\EssBdd{\hbox{\symbols n}\,}
\def\modSel{{\widetilde{J}}}

\par
\centerline{\bf Symmetry and short interval mean-squares}
\bigskip
\par
\centerline{\stampatello g. coppola - m. laporta}
\bigskip
\par
\noindent {\bf Abstract}. The weighted Selberg integral is a discrete mean-square, that is a generalization of the classical Selberg integral of primes to
	an arithmetic function $f$, whose values in a short interval are suitably attached to a weight function.
	We give conditions on $f$ and select a particular class of weights, in order to investigate 
	non-trivial bounds of weighted Selberg integrals of both $f$ and $f\ast\mu$. In particular, we discuss the cases of the symmetry integral and the modified Selberg integral, the latter  involving the Cesaro weight. We also prove some side results when $f$ is a divisor function.\smallskip

\bigskip

\par
\centerline{\stampatello 1. Introduction and statement of the results}
\smallskip
\par
\noindent
The {\it symmetry integral} of 
$f:\N\rightarrow\C$ is a short interval mean-square of the type
$$
J_{\sgn,f}(N,H)\defineq\avesum\Big|\sum_{x-H\le n\le x+H}f(n)\,\sgn(n-x)\Big|^2,
$$
\par
\noindent
where $N,H\in\N$ are such that $H=o(N)$ as $N\to \infty$, $x\sim N$ means that $x\in(N,2N]\cap\N$, and $\sgn(0)\defineq 0$, $\sgn(t)\defineq |t|/t$ for $t\not=0$. 
\par
These kinds of mean-squares have been intensively studied in [C1],[C2],[C4] and [C-L1] for different instances of $f$, with a particular attention to the divisor function $d_k$ for $k\ge 2$, where $d_k(n)$ is the number of ways to write $n$ as a product of $k$ positive integers. Symmetry integrals are so called because of the Kaczorowski and Perelli discover [K-P] about a remarkable link between the symmetry integral of the primes and the classical Selberg integral [Se], respectively 
$$
\int_{N}^{2N}\Big|\sum_{x-H\le n\le x+H}\Lambda(n)\sgn(n-x)\Big|^2 {\rm d}x, 
\enspace 
\int_{N}^{2N}\Big| \sum_{x<n\le x+H}\Lambda(n)-H\Big|^2 {\rm d}x,
$$
\par
\noindent
where $\Lambda$ is the von Mangoldt function, defined as $\Lambda(n)\defineq \log p$ if $n=p^r$ for some prime number $p$ and for some $r\in\N$, otherwise $\Lambda(n)\defineq 0$. In fact, 
even the discrete version of  the latter integral can be generalized 
to define the {\it Selberg integral} of any arithmetic function $f$, namely
$$
J_f(N,H)\defineq \avesum\Big|\sum_{x<n\le x+H}f(n)-M_f(x,H)\Big|^2, 
$$
\par
\noindent
where $M_f(x,H)$ is the so-called {\it short interval mean-value} of $f$.
\medskip
\par
\noindent
As in the prototype case of $d_k$, non-trivial estimates of $J_{\sgn,f}(N,H)$ are pursued for {\it essentially bounded} $f$, i.e. $f(n)\ll_{\varepsilon} n^{\varepsilon}\ \forall\varepsilon>0$. Plainly, the wider is the range of the {\it width} $\theta$ of the short interval $[x-H,x+H]$, with $H\asymp N^{\theta}$ (i.e $N^{\theta}\ll H\ll N^{\theta}$), for which a non-trivial bound holds, the finer is the result.
\medskip
\par
Our first theorem gives a link between non-trivial estimates
of the symmetry integrals of both $g$ and the Dirichlet convolution product $g\ast \1$, where $\1$ denotes the constantly $1$ function and $g$ is a real-valued and essentially bounded arithmetic function. Note that $g\ast \1$ is essentially bounded as well. 
\smallskip
\par
\noindent {\bf Theorem 1}.
{\it Let $g:\N\rightarrow\R$ be essentially bounded and $H\asymp N^{\vartheta}$ for a fixed $\vartheta \in (1/3,1)$. If there exists $G\in (0,1)$ such that
$$
J_{\sgn,g}(N,h)\ll Nh^{2-G}
$$
\par
\noindent
for every integer $h\asymp N^{\theta}$ with 
$$
\theta\in \left({{3\vartheta-1}\over {(1-G)\vartheta+G+1}},\vartheta\right],
$$
\par
\noindent
then				
$$
J_{\sgn,g\ast \1}(N,H)\ll NH^{2-G'}
$$
\par
\noindent
for every $\displaystyle{G'\in \left(0,\min\Big({{3-1/\vartheta}\over {1+2/G}},{1\over\vartheta}-1\Big)\right)}$.
}
\smallskip
\par
\noindent {\bf Remark 1}. {\it Positive {\sl exponent gains} $G$ and $G'$ convey  non-trivial bounds compared with $Nh^2$ and $NH^2$, respectively. For all $G\in (0,1)$ and $\vartheta \in (1/3,1)$ let us set 
$$
T_{\vartheta,G}\defineq  
\left({{3\vartheta-1}\over {(1-G)\vartheta+G+1}},\vartheta\right],
$$
\par
\noindent
and note that $T_{\vartheta,G}$ collapses to $\emptyset$, as $\vartheta\to 1^-$, while it enlarges to $(0,1/3]$, as $\vartheta\to 1/3^+$. Both limits hold
uniformly with respect to $G\in (0,1)$ and by Theorem 1 they yield $G'\to 0$, that is to say, both  $\vartheta\to 1^-$ and $\vartheta\to 1/3^+$ lead to a trivial bound for $J_{\sgn,g\ast \1}(N,H)$. Moreover, 
$$
\vartheta\in \left({1\over 3},{{G+1}\over{2G+1}}\right]
\Longleftrightarrow 
{{3-1/\vartheta}\over {1+2/G}}\le{1\over\vartheta}-1.
$$
\par
\noindent
In particular, for 
$$
\vartheta={{G+1}\over{2G+1}}\in\Big({2\over 3},1\Big)
$$
\par
\noindent
one gets the largest possible range for $G'\in\Big(0,\displaystyle{{{G}\over{G+1}}}\Big)$.
}

We postpone the proof of Theorem 1 until Sect.4 together with the proof of the following noteworthy consequence in the special cases of the functions $d_3$ and $\omega$, where
$\omega(n)\defineq\sum_{p|n}1$ counts the number of distinct prime divisors of $n$.
\smallskip
\par
\noindent {\bf Corollary 1}.
{\it For any integer $H\asymp N^{\vartheta}$, non-trivial bounds hold for both
$$
J_{\sgn,d_3}(N,H)\ \hbox{if} \enspace \vartheta \in \left({1\over 3},{1\over 2}-\varepsilon\right), \forall\varepsilon>0,
\quad 
J_{\sgn,\omega}(N,H)\ \hbox{if} \enspace \vartheta \in \left({7\over 17}+\varepsilon,1-\varepsilon\right), \forall\varepsilon>0.
$$
} 
\smallskip
\par
\noindent
After the introduction  of the {\it modified} Selberg integral [C0], 
$$
{\widetilde{J}}_f(N,H)\defineq \avesum\Big|{1\over H}\sum_{h\le H}\sum_{x-h<n<x+h}f(n)-M_f(x,H)\Big|^2, 
$$
\par
\noindent
where $M_f(x,H)$ is the same mean-value that appears in $J_f(N,H)$, we have further generalized mean-squares in short intervals for arithmetic functions gauged by a general weight $w$ (see [C-L]). Indeed, $J_{\sgn,f}(N,H)$, $J_f(N,H)$ and ${\widetilde{J}}_f(N,H)$ are particular instances of the so-called {\it weighted} Selberg integral of $f$, i.e.
$$
J_{w,f}(N,H)\defineq \avesum\Big|\sum_{n}\wH(n-x)f(n)-M_f(x,\wH)\Big|^2, 
$$
\par
\noindent
where $\wH$ is the product of the weight $w:\R\rightarrow \C$ and
the characteristic function $\oneH$ of the set $[-H,H]\cap\Z$, whereas
the mean value $M_f(x,\wH)$ has to be determined in agreement with the choice of  $w$. In particular, the weights involved in 
$J_f(N,H)$ and ${\widetilde{J}}_f(N,H)$ are respectively the {\it unit step} function
$$
u(a)\defineq
\cases{
1 & if $a>0$\cr 
0 & otherwise\cr} 
$$
\par
\noindent
and the {\it Cesaro weight} $\displaystyle{\CH(a)\defineq \max\Big(1-{|a|\over H},0\Big)}$.
While this is readily seen for $u$, the relation between ${\widetilde{J}}_f(N,H)$ and $\CH$
relies on
a well-known observation, due to the Italian mathematician Cesaro, i.e.
$$
\sum_{0\le |n-x|\le H}\left( 1-{{|n-x|}\over H}\right)f(n)={1\over H}\sum_{h\le H}\sum_{|n-x|<h}f(n).
$$
\par				
\noindent
Concerning the mean value terms, whereas it is plain that $M_f(x,\sgnH)$ vanishes identically for any $f$, according to Ivi\'c [Iv] if
$f$ has Dirichlet series $F(s)$ that is meromorphic in $\C$ and absolutely convergent in the half-plane $\Re(s)>1$  at least, then the mean value appearing in $J_f(N,H)$ and ${\widetilde{J}}_f(N,H)$  has the {\it analytic form} (see [C-L])
$$
M_f(x,H)= Hp_f(\log x),
$$
\par
\noindent
where $p_f(\log x)\defineq\Res_{s=1}F(s)x^{s-1}$ is the {\it logarithmic polynomial} of $f$.
More generally, in $J_{w,f}(N,H)$ one has (compare [C-L]) 
$$
M_f(x,\wH)=\FTwH(0)p_f(\log x),
$$
\par
\noindent
where
$\FTwH(0)\defineq \sum_a \wH(a)$ is the so-called
{\it mass} of $w$ in $[-H,H]$ (see Sect.2). From our study 
[C-L] it turns out that, if $f=g\ast \1$, then under suitable conditions on $g$ and $w$ one might  expect $M_f(x,\wH)$ 
to be close, at least in the mean-square, to its {\it arithmetic form} 
$$
\FTwH(0)\sum_{q\leq Q}{{g(q)}\over q},\ \hbox{for some}\  Q=Q(x)\ll N. 
$$
\par 
\noindent
This is a new feature exploited in
[C1],[C2] and [C-L1] when $f$ is a {\it sieve function}, that is
$$
f(n)\defineq \sum_{q|n\atop q\le Q}g(q), 
$$
\par
\noindent
where $g:\N \rightarrow \R$ is essentially bounded and $Q$ does not depend on $x\sim N$.
If $Q\ll N^{\lambda_f}$ for some $\lambda_f\in [0,1]$ and uniformly for $n\ll N$, then we say that $f$ has {\it level} (at most) $\lambda_f$. Trivially, any essentially bounded  arithmetic function is a sieve function of level at most $1$, while in the sequel a particular attention is given to arithmetic functions of level (strictly) less than $1$. 
A typical case where $Q$ depends on $x$ is the divisor function $d_k=d_{k-1}\ast \1$. The last section is devoted to such a case, 
in order to accomplish a discussion, commenced in [C-L], on the problem of showing sufficient proximity of the analytic and the arithmetic forms of the short interval mean value.\par
Here we are going to explore further the relation between weighted Selberg integrals of 
$g\ast \1$ and $g$. 
Our next result yields that (roughly speaking)
if $J_g(N,H)$ is close enough to ${\widetilde{J}}_g(N,H)$, then 
the same happens with $J_{g\ast\1}(N,H)$ and ${\widetilde{J}}_{g\ast\1}(N,H)$. 
\smallskip
\par
\noindent {\bf Theorem 2}.
{\it Let $g:\N\rightarrow\R$ be essentially bounded and let $w'$ be a weight such that $\wHone= \CH-\uH$ with
$H\asymp N^{\vartheta}$ for a fixed $\vartheta \in (1/3,1)$. If there exists $G\in (0,1)$ such that
$$
J_{w',g}(N,h)\ll Nh^{2-G}
$$
\par
\noindent
for every integer $h\asymp N^{\theta}$ with $\displaystyle{\theta\in \Big({{3\vartheta-1}\over {(1-G)\vartheta+G+1}},\vartheta\Big]}$,
then $$
J_{w',g\ast \1}(N,H)\ll NH^{2-3G'}\ \hbox{for some}\ G'=G'(\vartheta,G)>0. 
$$ 
} 
\smallskip
\par
\noindent
Since the proof goes as for Theorem 1 (in analogy with $\sgn$, the weight $w'$ has zero mass in $[-H,H]$), we omit it, while in Sect.4 we show that
Theorem 2 yields the aforementioned following consequence.
\smallskip
\par
\noindent {\bf Corollary 2}.
{\it Under the same hypotheses of Theorem 2 one has
$$
J_{g\ast \1}(N,H)-\modSel_{g\ast \1}(N,H)\ll NH^{2-G'}. 
$$
} 
\smallskip
\par				
The next theorem is a generalization to $J_{w,f}(N,H)$ of results in [C1] and [C-L1], when $f$ is a sieve function and $w$ belongs to a particular class of weights, that we describe as follows. If the correlation of $\wH$, i.e.
$$
\Corr_{\wH}(a)\defineq \doublesum_{{m \thinspace \quad \thinspace n}\atop {m-n=a}}\wH(m)\overline{\wH(n)}, 
$$
\par
\noindent
satisfies the formula
$$
\sum_{a\equiv 0\, (\!\!\bmod \thinspace \ell)}\Corr_{\wH}(a)={1\over {\ell}}\sum_{a}\Corr_{\wH}(a)+O(H)
\quad \forall \ell\le 2H,
\leqno{(A)}
$$
\par
\noindent
then $\wH$ is said to be {\it arithmetic} and, if this is the case for every $H$, then
we say so for $w$. In particular, $w$ is a {\it good}
weight if it is arithmetic and
{\it absolutely bounded} on the integers, i.e. there exists $K\in(0,+\infty)$ such that
$|w(n)|\le K \enspace \forall n\in \Z$.
In Sect.2 we show that the weights $u,\sgn$ and $\CH$ are 
arithmetic (good weights {\it de facto}, it being plain that they are absolutely bounded by $1$).
\smallskip
\par
In order to state the next theorem, we need also some further notation and convention. 
First, the following modified version of Vinogradov's notation is useful to hide arbitrarily small powers (of the main variable):
$$
A\EssBdd B
\enspace \defin \enspace 
A\ll_{\varepsilon} N^{\varepsilon}B,
\quad \forall \varepsilon>0
$$
\par
\noindent
Then, by writing $\supporto g\subseteq [1,Q]$ we implicitly mean that the arithmetic function under consideration is $g\cdot \1_{[1,Q]}$, where
$\1_{[1,Q]}$ is the characteristic function of $[1,Q]\cap\N$, so that $g\ast \1$ is a sieve function.
\smallskip
\par
\noindent {\bf Theorem 3}.
{\it Let $g:\N\rightarrow\R$ be essentially bounded and such that  $\supporto g\subseteq [1,Q]$ with
$Q=Q(N,H)\to \infty$. For every good and real weight $w$ one has
$$
J_{w,g\ast \1}(N,H)\EssBdd NH+Q^2H+QH^2+H^3.
$$
} 
\smallskip
\par
\noindent
The proof is postponed until Sect.4. Here we point out that, under the same hypotheses for $g$ and $w$, if 
$H\asymp N^{\vartheta}$ with $\vartheta\in(0,1/2)$, then we have recently 
established (see [C-L1])
\smallskip
\par
\centerline{$J_{w,g\ast \1}(N,H)\EssBdd NH+N^{\delta}Q^{95/48}H^2+N^{1-2\delta/3}H^2+QH^2 \qquad \forall\delta>0$,}
\smallskip
\par
\noindent
by means of a very technical result based upon averages of Kloosterman sums [D-F-I]. 
\bigskip

The last result of this section exhibits a {\it length-inertia} property for the Selberg integral, that allows to preserve non-trivial bounds, as the length of the short interval increases. 
\par
Let $[x]$ be the integer part of $x\in\R$ and let $L\defineq \log N$. In Sect.4 we prove also the following theorem.
\smallskip
\par
\noindent {\bf Theorem 4}.
{\it Let $f:\N\rightarrow\R$ be such that the logarithmic polynomial $p_f(\log n)$ is defined. For every $H>h$ one has
$$
J_f(N,H)\ll H^2h^{-2} J_f(N,h) + J_f(N,H-h[H/h]) + H^3 \left( \Vert f\Vert_{\infty}^2 + L^{2c}\right),
$$
\smallskip
\par
\noindent
where 
$$
\Vert f\Vert_{\infty}\defineq \max_{[N-H,2N+H]}|f|
$$
\par
\noindent
and $c$ is the degree of $p_f$.
} 
\smallskip
It is worthwhile to remark that, while the symmetry integral has no length-inertia property, we have recently showed  that
such a feature holds also for the modified Selberg integral. Indeed, in [C-L2] by means of our version of a Gallagher's inequality we have proved that
$$
\modSel_f(N,H)\ll H^2h^{-2} \modSel_f(N,h) 
              + \left(Nh^4H^{-2} 
               + H^3\right)  \Vert f\Vert_{\infty}^2 + H^3L^{2c},
$$
\par
\noindent
under the same hypotheses of Theorem 4. We underline that no Gallagher type inequality is used in the proof of Theorem 4, that is proved in elementary fashion. 
\bigskip
\par				
\noindent
{\bf Plan of the paper.} In Sect.2 we introduce the {\it sporadic} functions and discuss some properties of the arithmetic weights.
The necessary Lemmata for Theorems 1 and 3 constitute the third section of the paper, whereas Sect.4 contains the proofs of Theorems 1, 3, 4 and Corollaries 1, 2. The last section complements our study in [C-L] dealing with the case of the divisor function $d_k$. 
\bigskip
\par
\noindent
{\bf Notation and definitions.} If the implicit constants in the symbols $O$ and $\ll$ depend on 
some parameters like $\varepsilon>0$, then mostly we specify it by introducing subscripts
like $O_\varepsilon$ and $\ll_\varepsilon$, whereas we omit subscripts 
for $\EssBdd$ defined above. 
Notice that the value of $\varepsilon$ may change from statement to
statement, since $\varepsilon>0$ is arbitrarily small. For the sake of clarity, let us remark that throughout the paper $H,N$ denote positive integers such that $H=o(N)$ as $N\to \infty$, i.e.
$H/N\to 0$ as $N\to \infty$. The notation $H=o(N)$ is used synonymously with $N=\infty(H)$.
Typically we write
$H\asymp N^{\theta}$ for $N^{\theta}\ll H\ll N^{\theta}$ with
$\theta\in(0,1)$, that we call the {\it width} of the short interval $[x-H,x+H]$.
As already mentioned, we use to abbreviate $L=\log N$. 
\par
The symbol $\1$ denotes the constantly $1$ function, while $\1_{U}$ is the characteristic function of $U\cap\Z$ for every $U\subseteq\R$.
In particular, we abbreviate $\oneH=\1_{[-H,H]}$, so that $\wH=w\cdot \oneH$ for $w:\R\rightarrow\C$.  
If
$f=g\ast \1$, then $g=f\ast\mu$ is called the {\it Eratosthenes transform} of $f$. Therefore, $\1$  is the Eratosthenes transform of the divisor function $\divisor=\1 \ast \1$. 
More generally, $d_k=\underbrace{\1 \ast \cdots \ast \1}_{k\, {\rm times}}=d_{k-1}\ast\1$ for $k\ge 3$. 
\par
\noindent
In sums like $\sum_{a\le X}$ it is implicit that $a\ge 1$, while the range of the $a$'s 
in $\sum_{a}$ is the support of the function that appears in the summands. As usual, $e(\alpha)\defineq e^{2\pi i\alpha},\ \forall \alpha \in \R$, and 
$e_q(a)\defineq e(a/q),\ \forall (q,a)\in \N\times\Z$. Especially within formulae, at times we abbreviate $n\equiv 0\ (\bmod\, q)$ as
$n\equiv 0\ (q)$. The symbol $\starsum$ indicates that the sum
is taken over the reduced residues. 
The distance of $\alpha \in \R$ from the nearest integer is 
$\Vert\alpha\Vert\defineq \min\big(\{ \alpha\}, 1-\{ \alpha\}\big)$, where $\{ \alpha\}\defineq \alpha-[\alpha]$ is the fractional part of $\alpha$. 
Throughout the paper we apply standard formulae without further references. 
For example, we use the asymptotic equation ($\gamma$ is the {\it Euler-Mascheroni constant})
$$
\sum_{n\le x}{1\over n}=\log x+\gamma+O(1/x).
$$

\bigskip

\par
\centerline{\stampatello 2. Sporadic functions and arithmetic weights}
\smallskip
\par
\noindent
As already mentioned, any function of the type $w:\R \rightarrow \C$ can take over the role of a weight here, though {\it de facto}
we deal with the weighted characteristic function $\wH=w\cdot\oneH$ of the integers in the short interval $[-H,H]$. We consider the associated exponential sum
$$
\FTwH(\beta)\defineq \sum_{a}\wH(a)e(a\beta)
\quad 
\forall \beta \in [0,1),
$$
\par
\noindent
whose value for $\beta=0$ is the {\it mass} of $w$ in $[-H,H]$, i.e. 
$$
\FTwH(0)=\sum_a \wH(a).
$$
\par
\noindent
Since $\#\{n\in[x-H,x+H]:\ n\equiv 0\ (\bmod \, q)\}\le 1$ when $q>2H$, then we call
$$
{\cal W}_H(x;q)\defineq \sum_{n\equiv 0\, (q)}\wH(n-x)=\sum_{a\equiv -x\, (q)}\wH(a), 
$$
\par
\noindent
the (weighted) {\it sporadic sum}, while the (weighted) {\it sporadic function} is 
$$
\chi_q(x,\wH)\defineq {\cal W}_H(x;q)-{\FTwH(0)\over q}= \sum_{a\equiv -x\, (q)}\wH(a)-{1\over q}\sum_{a} \wH(a),
$$
\par
\noindent
where $\FTwH(0)/q$ exhibits a behavior of a {\it local} mean-value for ${\cal W}_H(x;q)$.
\smallskip
\par				
\noindent
In the next proposition we give the so-called {\it Fourier-Ramanujan expansion} of 
the sporadic function. 
\smallskip
\par				
\noindent {\bf Proposition 1}.
{\it For every $w:\R \rightarrow \C$ and all $q,H\in\N$ one has
$$
\chi_q(x,\wH)={1\over q}\sum_{d>1\atop d|q}\starsum_{j\le d}\FTwH\Big( {j\over d}\Big)e_{d}(jx)
$$
\par
\noindent
(we assume that the sum vanishes for $q=1$).
} 
\smallskip
\par
\noindent {\bf Proof}. From the orthogonality of the additive characters, 
$$
{1\over q}\sum_{r\le q}e_q(ar)=
\cases{1 & if $q|a$,\cr 
0 & otherwise\ .
\cr
}
$$
\par
\noindent
one gets
$$
\sum_{n\equiv 0\, (q)}\wH(n-x)=\sum_{a+x\equiv 0\, (q)}\wH(a)
={1\over q}\sum_{r\le q} e_q(rx)\sum_{a}\wH(a)e_q(ar)
={1\over q}\sum_{r\le q}\FTwH\left({r\over q}\right)e_q(rx).
$$
\par
\noindent
Thus, since $\FTwH(0)=\FTwH(1)$, we can write
$$
\chi_q(x,\wH)=\sum_{n\equiv 0\, (q)}\wH(n-x)-{\FTwH(0)\over q}
={1\over q}\sum_{r<q}\FTwH\left( {r\over q}\right)e_q(rx)
={1\over q}\sum_{d>1\atop d|q}\sum_{{r<q}\atop {(r,q)=q/d}}\FTwH\left( {r\over q}\right)e_q(rx) 
$$
$$
={1\over q}\sum_{d>1\atop d|q}\sum_{{j\le d}\atop {(j,d)=1}}\FTwH\left( {j\over d}\right)e_d(jx).
$$
\par
\noindent 
The proposition is proved.\hfill $\square$ 
\smallskip
\par
Recall that $\wH$ is {\it arithmetic} if the correlation $\Corr_{\wH}(a)$ satisfies $(A)$. 
The next propositions give some properties of such weights.
\smallskip
\par
\noindent {\bf Proposition 2}.
{\it For every $w:\R \rightarrow \C$  and every $H\in\N$ the following properties hold.
\item{(i)} $\wH$ is arithmetic if and only if
$$
{1\over {q}}\sum_{j<q}\Big|\FTwH\left({j\over q}\right)\Big|^2\ll H
\quad \forall q\le 2H.
\leqno{(B)}
$$
\item{(ii)} If $w$ is a good weight, then $(B)$ holds for all $q\ge 1$.
\item{(iii)} If $\wH$ is arithmetic, so is its normalized correlation ${\Corr_{\wH}/H}$.
} 
\smallskip
\par
\noindent {\bf Proof}. (i) Similarly to Proposition 1, through an application of the orthogonality of additive characters it easily seen that $(A)$ is equivalent to 
$$
{1\over {q}}\sum_{j<q}\widehat{\Corr}_{\wH}\Big({j\over q}\Big)\ll H
\quad \forall q\le 2H.
$$
\par
\noindent
Therefore, $(B)$ follows immediately from the identity
$$
\widehat{\Corr}_{\wH}(\beta)=\sum_a \Corr_{\wH}(a)e(a\beta)
=\sum_a \doublesum_{m-n=a}\wH(m)\overline{\wH(n)}e(a\beta)
=\Big|\sum_r \wH(r)e(r\beta)\Big|^2=|\FTwH(\beta)|^2\
\forall \beta \in [0,1).
$$
\par
\noindent 
(ii) Note that if $w$ is absolutely bounded, then $(B)$ holds for {\it high divisors} $q>2H$ because
$$
{1\over {q}}\sum_{j<q}\Big| \FTwH\Big({j\over q}\Big)\Big|^2\le {1\over {q}}\sum_{j\le q}\Big| \FTwH\Big({j\over q}\Big)\Big|^2
={1\over {q}}\sum_{j\le q}\sum_{h_1}\sum_{h_2}\wH(h_1)\overline{\wH(h_2)}e_q(j(h_1-h_2))= 
$$
$$				
=\sum_{0\le |h_1|\le H}\sum_{{0\le |h_2|\le H}\atop {h_2\equiv h_1\, (q)}}\wH(h_1)\overline{\wH(h_2)}
=\sum_{0\le |h_1|\le H}w(h_1)\sum_{{0\le |h_2|\le H}\atop {h_2\equiv h_1\, (q)}}\overline{w(h_2)}\ll {{H^2}\over q}+H.  
$$
\par
\noindent
(iii) Let us show that $(B)$ holds for $\Corr_{\wH}/H$. Indeed, since trivially 
$\FTwH(\beta)\ll H$, 
$\forall \beta\in[0,1)$, 
 then
$$
{1\over {q}}\sum_{j<q}\Big|\widehat{\Big({\Corr_{\wH}\over H}\Big)}\Big({j\over q}\Big)\Big|^2={1\over {qH^2}}\sum_{j<q}\Big|\widehat{\Corr}_{\wH}\Big({j\over q}\Big)\Big|^2
={1\over {qH^2}}\sum_{j<q}\Big|\FTwH\Big({j\over q}\Big)\Big|^4
\ll {1\over {q}}\sum_{j<q}\Big|\FTwH\Big({j\over q}\Big)\Big|^2
\ll H
\quad \forall q\le 2H.
$$
\par
\noindent 
The proposition is completely proved.\hfill $\square$ 
\smallskip
\par
\noindent 
Now, we show that $(B)$ holds for $u,\sgn$ and $\CH$, that is to say, such weights are good.
To this end, we set
$$
{\cal L}^{2}_q(\FTwH)\defineq 
{1\over {q^2}}\sum_{j<q}\Big| \FTwH\left({j\over q}\right)\Big|^2
$$
\par
\noindent 
and prove the next property that plainly implies $(B)$ for $\uH, \sgnH$ and $\CH$.
\smallskip
\par
\noindent {\bf Proposition 3}.
{\it If $\wH\in\{\uH, \sgnH, \CH\}$, then} 
$$
{\cal L}^{2}_q(\FTwH)\ll \min(1,H/q).
$$
\smallskip
\par
\noindent {\bf Proof}. Let us start with $\wH=\uH$, whose associated exponential sum satisfies the well known inequality 
$$
|\FTuH(\alpha)|=\Big|\sum_{n\le H}e(n\alpha)\Big|=\Big|{{\sin(\pi H\alpha)}\over {\sin(\pi \alpha)}}\Big|\ll 
\min\Big(H,{1\over {\Vert \alpha \Vert}}\Big).
$$
\par
\noindent 
We set $H_q\defineq q\{ H/q\}\ll \min(q,H)$ and note that $\{ H/q\}=0$ would yield that $\sin(\pi Hj/q)=0$ and the stated inequality would hold trivially for ${\cal L}^{2}_q(\FTuH)=0$. Thus, we can assume $H_q>0$ to see that 
$$
{\cal L}^2_q(\FTuH)={1\over {q^2}}\sum_{j<q}\Big|{{\sin(\pi jH_q/q)}\over {\sin(\pi j/q)}}\Big|^2
\ll {1\over {q^2}}\sum_{0<|j|\le {q\over {H_q}}}H_q^2 + \sum_{{q\over {H_q}}<|j|\le {q\over 2}}{1\over {j^2}} 
\ll {{H_q}\over q} 
\ll \min\Big(1,{H\over q}\Big).
$$
\par
\noindent 
Since the Cesaro weight is the normalized correlation of $u$, i.e. 
$$
\CH(a)
={1\over H} \sum_{t\le H-|a|}1 ={1\over H} \doublesum_{{m,n\le H}\atop {m-n=a}}1
={\Corr_{\uH}(a)\over H},
$$ 
\par
\noindent 
then $\CH$ is arithmetic because of the third assertion in Proposition 2. However, it turns out that the stated inequality holds also in this case:
$$
{\cal L}_q^2(\FTCH)=
{\cal L}^2_q\Big(\widehat{{{\Corr_{\uH}}\over H}}\Big)={1\over {H^2}}{\cal L}^2_q(|\FTuH|^2)
\ll {\cal L}^2_q(\FTuH)
\ll \min\Big(1,{H\over q}\Big).
$$
\par
\noindent
Now, let us consider the case $\wH=\sgnH$ and recall that for every $\alpha \in\R\setminus\Z$ we can write (see [Da, Ch.25]) 
$$
\left|\FTsgnH(\alpha)\right|=2\Big|\sum_{h\le H}\sin(2\pi h\alpha)\Big|
\ll \left|\cot(\pi \alpha)\right| \sin^2(\pi H\alpha) + \left|\sin(2\pi H\alpha)\right|
\ll {{\sin^2(\pi H\alpha)}\over {|\sin(\pi \alpha)|}}
\ll {{\Vert H\alpha\Vert^2}\over {\Vert \alpha\Vert}}, 
$$
\par
\noindent 
while $\FTsgnH(\alpha)=0$ when $\alpha \in \Z$. Thus, as above one gets (compare [C-S]) 
$$
{\cal L}^2_q(\FTsgnH)\ll \sum_{0<|j|\le q/2} {{\Vert Hj/q\Vert^4}\over {|j|^2}}
\ll \Big({H_q\over {q}}\Big)^4\sum_{0<|j|\le {q\over {2H_q}}}j^2 + \sum_{{q\over {2H_q}}<|j|\le {q\over 2}}{1\over {j^2}}
\ll {{H_q}\over q}
\ll \min\Big(1,{H\over q}\Big). 
$$
\par				
\noindent
The proposition is completely proved.\hfill $\square$ 
\smallskip
\par
\noindent {\bf Remark 2}. {\it Similarly to the unit step function $u$, it is not difficult to see that any piecewise-constant weight is arithmetic.
Moreover, as for $u$, it is plain that for the normalized correlation of $\sgn$ we get 
$$
{\cal L}^2_q\Big(\widehat{{{\Corr_{\sgnH}}\over H}}\Big)={1\over {H^2}}{\cal L}^2_q(|\FTsgnH|^2)
\ll \min\Big(1,{H\over q}\Big). 
$$
}
\par
\noindent
We close this section with an example of an absolutely bounded weight that is not arithmetic. To this end, let us take $\wH(m)=e(m\eta)\uH(m)$ for a fixed $\eta\in[0,1)$ to be chosen later. Then, its correlation is 
$$
\Corr_{\wH}(a)=\sum_{{0<m\le H}\atop {0<m-a\le H}}e(m\eta)e(-(m-a)\eta)
=e(a\eta)\sum_{m_1<m\le m_2}1=e(a\eta)\max(H-|a|,0), 
$$
\par
\noindent
where $m_1\defineq\max(0,a), m_2\defineq\min(H,H+a)$.  Note that $\Corr_{\wH}$
satisfies the formula
$$
\sum_{a\equiv 0\, (\ell)}\Corr_{\wH}(a)=\sum_{0\le |b|\le [H/\ell]}(H-\ell |b|)e(\ell b\eta)
=\ell \sum_{0\le |b|\le [H/\ell]}([H/\ell]-|b|)e(\ell b\eta)+O(H)= 
$$
$$
=\ell \Big|\sum_{h\le [H/\ell]}e(h\ell \eta)\Big|^2+O(H). 
$$
\par
\noindent
In particular, by taking $\eta=m/\ell$ with $1\le m<\ell$ one has
$$
\sum_{a\equiv 0\, (\ell)}\Corr_{\wH}(a)=\ell \Big|\sum_{h\le [H/\ell]}e(hm)\Big|^2+O(H)=
{{H^2}\over {\ell}}+O(H), 
$$
\par
\noindent 
so that
$$
\sum_{a\equiv 0\, (\ell)}\Corr_{\wH}(a)=\infty(H)\Longleftrightarrow \ell=o(H). 
$$
\par
\noindent
Since (see the proof of (i) in Proposition 2)
$$
\sum_a \Corr_{\wH}(a)=\widehat{\Corr}_{\wH}(0)=|\FTwH(0)|^2=\Big|\sum_{a\le H}e(a\eta)\Big|^2, 
$$
\par
\noindent
then for $\eta=m/\ell$ we get
$$
\sum_a \Corr_{\wH}(a)=\Big|\sum_{a\le H}e_{\ell}(ma)\Big|^2
=\Big|\sum_{a\le \ell \{ H/\ell \} }e_{\ell}(ma)\Big|^2\ll \ell^2. 
$$
\par
\noindent
Hence, $(A)$ cannot hold for any choice of $\ell=o(H)$.

\bigskip

\par
\centerline{\stampatello 3. Lemmata for Theorems 1 and 3}
\smallskip
\par
\noindent
Here we prove the necessary lemmata for Theorems 1 and 3. More precisely, the first three
lemmas are applied within the proof of Theorem 1, while Lemma 4 is of use for Theorem 3. 
\smallskip
\par
\noindent {\bf Lemma 1}.
{\it Let $\kappa:x\in (N,2N]\rightarrow \kappa(x)\in[0,+\infty)$ be strictly increasing and such that $\kappa(2N)\ll Q\ll N$, where $Q$ may depend on $N, H$.
For every essentially bounded $g:\N \rightarrow \C$ and
every good weight $w$
one has 
$$
\avesum \Big|\sum_{q\le \kappa(x)}g(q)\chi_q(x,\wH)\Big|^2\EssBdd (N+Q^2)H\, .
$$
} 
\smallskip
\par				
\noindent {\bf Proof}. First, let us introduce the following auxiliary notation:
$$
{\cal W}^\ast_H(x;q)\defineq {1\over q}
\starsum_{j\le q}\FTwH
\Big({j\over q}\Big)e_q(jx)
$$
\par
\noindent
Then, by applying Proposition 1 we see that the left hand side of the inequality to be proved is equal to
$$
\doublesum_{d,q\le \kappa(2N)}\overline{g(d)}g(q)\sum_{{x\sim N}\atop {{x\ge \kappa^{-1}(d)}\atop {x\ge \kappa^{-1}(q)}}}
{1\over q}\sum_{\ell|q,\ell>1}\ell{\cal W}^\ast_H(x;\ell)
{1\over d}\sum_{t|d,t>1}t\overline{{\cal W}^\ast_H(x;t)}= 
$$
$$
=\doublesum_{1<\ell,t\le \kappa(2N)} {\cal W}^\ast_H(0;\ell)
   \overline{{\cal W}^\ast_H(0;t)} 
\sum_{n\le \kappa(2N)/\ell}{{\overline{g(\ell n)}}\over n}\sum_{m\le \kappa(2N)/t}{{g(tm)}\over m}
\sum_{{x\sim N}\atop {{x\ge \kappa^{-1}(\ell n)}\atop {x\ge \kappa^{-1}(tm)}}}e(\delta x), 
$$
\par
\noindent
where we set 
$\displaystyle{
\delta \defineq \Big\Vert {j\over {\ell}}-{r\over t}\Big\Vert \in \Big[0,{1\over 2}\Big]
}$.              
Because of the coprimality conditions, the {\it diagonal} terms are characterized by
$\delta=0\Leftrightarrow (j=r, \ell=t)$,
while $\delta \neq 0$ implies the classical well-spaced condition on the {\it Farey fractions} $j/\ell, r/t$,
$$
\delta = \Big\Vert {{jt-r\ell}\over {\ell t}}\Big\Vert \ge {1\over {\ell t}} \ge {1\over {\kappa(2N)^2}}.
$$
\par
\noindent
Now, the hypothesis that $w$ is a good weight yields
${\cal L}^{2}_q(\FTwH)\ll H/q\ \forall q\ge 1$ (see Proposition 2), which in turn implies that
$$
{1\over {q^2}}\starsum_{j\le q}\Big| \FTwH\Big({j\over q}\Big)\Big|^2
\ll {H\over q}\quad \forall q\ge 1. 
$$
\par
\noindent
Thus, 
by using the well-known estimate [Da, Ch.25] 
$$
\sum_{{x\sim N}\atop {{x\ge \kappa^{-1}(\ell n)}\atop {x\ge \kappa^{-1}(tm)}}}e(\delta x)\ll \min\Big(N, {1\over {\delta}}\Big) 
$$
\par
\noindent
and by applying the Large Sieve inequality in the form given in [C-S, Lemma 2], we conclude that
$$
\avesum \Big|\sum_{q\le \kappa(x)}g(q)\chi_q(x,w)\Big|^2 
\EssBdd (N+Q^2)\sum_{1< d \ll Q}{1\over {d^2}}\starsum_{j\le d}\Big|\FTwH \Big({j\over {d}}\Big)\Big|^2
\EssBdd (N+Q^2)H\sum_{1< d \ll Q}{1\over {d}} 
\EssBdd (N+Q^2)H.
$$
\par
\noindent
The lemma is proved.\hfill $\square$ 
\smallskip
\par
\noindent {\bf Remark 3}. {\it We explicitly note that under the same hypothesis for $g$ and $w$ through a similar proof one gets 
$$
\avesum \Big|\sum_{q\sim Q}g(q)\chi_q(x,\wH)\Big|^2\EssBdd (N+Q^2)H\quad\forall Q\ll N. 
$$
\par
\noindent
It is this inequality that is applied within the proof of Theorem 1.
}
\smallskip
\par
\noindent {\bf Lemma 2}.
{\it Let $A, B, Q, H, N$ be positive real numbers such that, as $N\to \infty$, $H=o(N)$, $H\to \infty$, and 
$$
Q\ll A<B\ll Q\ll N, \hbox{\rm with} \enspace Q=\infty(N/H). 
$$
\par
\noindent
For every essentially bounded $g:\N \rightarrow \C$ and every absolutely bounded weight $w$ with
$\FTwH(0)=0$, 
one has 
$$
\avesum\Big|\sum_{A<q\le B}g(q)\chi_q(x,\wH)\Big|^2
\EssBdd 
\avesum\Big|\sum_{{x\over B}<m\le {x\over A}}\sum_{q}g(q)\wH(mq-x)\Big|^2 
+ H^3.
$$
} 
\par				
\noindent {\bf Proof}. First, notice that $\FTwH(0)=0$ implies 
$$
\chi_q(x,\wH)={\cal W}_H(x;q)=\sum_{n\equiv 0\, (q)}\wH(n-x). 
$$ 
\par
\noindent
Then, let us apply Dirichlet's hyperbola method to get 
$$
\sum_{A<q\le B}g(q)\chi_q(x,\wH)=\sum_{A<q\le B}g(q)\sum_{n\equiv 0\, (q)}\wH(n-x)
=\sum_{{{x-H}\over B}\le m<{{x+H}\over A}}\sum_{{A<q\le B}\atop {{x-H}\over m}\le q\le {{x+H}\over m}}g(q)w(mq-x)= 
$$
$$
=\sum_{{{x+H}\over B}\le m<{{x-H}\over A}}\sum_{{{x-H}\over m}\le q\le {{x+H}\over m}}g(q)w(mq-x) + O_{\varepsilon}\Big(N^{\varepsilon}\sum_{m\in{\cal H}}
{\cal I}_H(x;m)\Big)= 
$$
$$ 
=\sum_{{x\over B}<m\le {x\over A}}\sum_{q}g(q)\wH(mq-x) + O_{\varepsilon}\Big(N^{\varepsilon}\sum_{m\in{\cal H}}{\cal I}_H(x;m)\Big), 
$$
\par
\noindent
where we have set 
$$
{\cal H}={\cal H}(x,H,A,B)\defineq \Big(\big[{{x-H}\over B},{{x+H}\over B}\big)\cup\big[{{x-H}\over A},{{x+H}\over A}\big)\Big)\cap\N,
$$
\par
\noindent
and it turns out that
$$
{\cal I}_H(x;m)\defineq\sum_{a\equiv -x\, (m)}\oneH(a)=\sum_{{{x-H}\over m}\le q\le {{x+H}\over m}}1\leq {{2H}\over m}+1\ll {H\over m}
$$
\par
\noindent 
from the hypothesis  $Q=\infty(N/H)$. Thus, the lemma is proved once we show that the $m-$sum in the $O_{\varepsilon}$-term gives a mean-square contribution $\ll H^3$. To this end, we note that $H/m \ll QH/N$ and get 
$$
\avesum \Big|\sum_{m\in{\cal H}}{\cal I}_H(x;m)\Big|^2 
\ll \Big({{QH}\over N}\Big)^2 \avesum \Big(\Big|\sum_{{{x-H}\over A}\le m\le {{x+H}\over A}}1\Big|^2 
                                           + \Big|\sum_{{{x-H}\over B}\le m\le {{x+H}\over B}}1\Big|^2\Big)\ll 
$$
$$
\ll \Big({{QH}\over N}\Big)^2 \Big(\sum_{{{N-H}\over A}\le m_1\le {{2N+H}\over A}}\sum_{m_1-{{2H}\over A}\le m_2\le m_1+{{2H}\over A}}1 + 
    + \sum_{{{N-H}\over B}\le m_1\le {{2N+H}\over B}}\sum_{m_1-{{2H}\over B}\le m_2\le m_1+{{2H}\over B}}1\Big)H\ll 
$$
$$
\ll \Big({{QH}\over N}\Big)^2 {N\over Q}\Big({H\over Q}+1\Big)H\ll \Big({H\over N}+{Q\over N}\Big)H^3\ll H^3. 
$$
\par
\noindent
The lemma is completely proved.\hfill $\square$ 
\smallskip
\par
\noindent {\bf Lemma 3}.
{\it Let $A,B,Q,H,N$ be as in Lemma 2. For every  essentially bounded $g:\N \rightarrow \C$, we get 
$$
\avesum \Big|\sum_{{x\over B}<m\le {x\over A}}\sum_{q}g(q)\sgnH(mq-x)\Big|^2 
\EssBdd {N\over Q} \sum_{m\asymp {N\over Q}}m\sum_{n\sim{{N}\over m}}\Big|\sum_{q}g(q)\sgn_{\big[{H\over m}\big]}(q-n)\Big|^2
         + {{N^3}\over {Q^2}} + H^3.
$$
} 
\smallskip
\par
\noindent {\bf Proof}. By the Cauchy inequality we see that the left hand side of the inequality to be proved is 
$$
\sum_{{N\over B}<m_1,m_2\le {{2N}\over A}}
\sum_{{x\sim N}\atop {{m_1A\le x<m_1B}\atop {m_2A\le x<m_2B}}}
\sum_{q_1}g(q_1)\sgnH(m_1q_1-x)
\sum_{q_2}\overline{g(q_2)\sgnH(m_2q_2-x)}\ll 
$$
$$
\ll \sum_{m_1,m_2\asymp {N\over Q}}\sqrt{\avesum\Big|\sum_{q_1}g(q_1)\sgnH(m_1q_1-x)\Big|^2 \avesum\Big|\sum_{q_2}g(q_2)\sgnH(m_2q_2-x)\Big|^2}\ll 
$$
$$				
\ll {N\over Q}\sum_{m\asymp {N\over Q}}\avesum\Big|\sum_{q}g(q)\sgnH(mq-x)\Big|^2. 
$$
\par
\noindent
Since
$$
\sum_{q}g(q)\sgnH(mq-x)=\sum_{{{x-H}\over m}\le q\le {{x+H}\over m}}g(q)\sgn(mq-x)
=\sum_{q}g(q)\sgn_{[H/m]}(q-[x/m])+O_{\varepsilon}(N^{\varepsilon}), 
$$
\par
\noindent
we write
$$
\avesum \Big|\sum_{{x\over B}<m\le {x\over A}}\sum_{q}g(q)\sgnH(mq-x)\Big|^2
\EssBdd {N\over Q}\sum_{m\asymp {N\over Q}}\avesum \Big|\sum_{q}g(q)\sgn_{\big[{H\over m}\big]}\big(q-\big[{x\over m}\big]\big)\Big|^2+{{N^3}\over {Q^2}}.
$$
\par
\noindent
Thus, by writing $x=mn+r$ with $0\le r\le m-1$ (compare [C5, Lemma 2.4]), we see that
$$
\avesum\Big|\sum_{q}g(q)\sgn_{[H/m]}(q-[x/m])\Big|^2 \EssBdd m\sum_{n\sim {{N}\over m}}\Big|\sum_{q}g(q)\sgn_{[H/m]}(q-n)\Big|^2 + m[H/m]^2,
$$
\par
\noindent
where $m[H/m]^2\ll H^2/m$ gives clearly a contribution $\EssBdd H^3$.\hfill $\square$ 
\smallskip
\par
\noindent {\bf Lemma 4}.
{\it Let $g:\N \rightarrow \R$ be essentially bounded and such that $\supporto g\subseteq [1,Q]$ with $Q=Q(N,H)\ll N$
and let $f=g\ast \1$.
\par
\noindent
(i) For every 
absolutely bounded weight $w:\R\rightarrow\R$ one has
$$
J_{w,f}(N,H)-\sum_{0\le |a|\le 2H}\Corr_{\wH}(a)\Corr_f(a)+
N\FTwH(0)^2\Big(\sum_{q\le Q}{{g(q)}\over q}\Big)^2
\EssBdd H^3+H^2Q. 
$$
\par
\noindent
(ii) For every integer $a\neq 0$ one has
$$
\Corr_f(a)=\sum_{\ell|a}\doublesum_{(d,q)=1}g(\ell d){{g(\ell q)}\over q}\left( \left[ {{2N}\over {\ell d}}\right]-\left[ {N\over {\ell d}}\right]\right)+R_f(a)
$$
\par
\noindent
with
$$
R_f(a)\defineq \sum_{\ell|a}\doublesum_{(d,q)=1}g(\ell d){{g(\ell q)}\over q}
\dashsum_{j(q)}e_q\Big(-{ja\over\ell}\Big)\sum_{m\sim {N\over {\ell d}}}e_q(jdm), 
$$
\par
\noindent
where the dashed sum is over all the nonzero residue classes $j$ {\rm mod} $q$.
\par
\noindent
Further, if $K:[-2H,2H]\rightarrow \C$ is an even function, then
$$
\sum_{a\neq 0} K(a)R_f(a)=\sum_{\ell \le 2H}\doublesum_{(d,q)=1}g(\ell d){{g(\ell q)}\over q}\dashsum_{j(q)}\sum_{m\sim {N\over {\ell d}}}\cos{{2\pi jdm}\over q}\sum_{a\neq 0}K(a\ell)e_q(ja). 
$$
\par
\noindent
(iii) For every absolutely bounded weight $w:\R\rightarrow\R$ one has
$$
\sum_{a\neq 0}\Corr_{\wH}(a)R_f(a)\EssBdd NH + Q^2 H + QH^2. 
$$
} 
\smallskip
\par
\noindent {\bf Proof}. (i) is a consequence of Lemma 7 in [C-L], that precisely yields the formula
$$
J_{w,f}(N,H)=\sum_{0\le |a|\le 2H}\!\!\Corr_{\wH}(a)\Corr_f(a)-2\sum_n f(n)\avesum \wH(n-x)M_f(x,\wH)
             +\avesum M_f(x,\wH)^2+O\left(H^3\Vert f\Vert_{\infty}^2\right), 
$$
\par				
\noindent
where $\displaystyle{\Vert f\Vert_{\infty}\defineq \max_{N-H<n\le 2N+H}|f(n)|}$.
Indeed, since for the mean-value of the sieve function $f$ one has 
$$
M_f(x,\wH)=\FTwH(0)\sum_{q\le Q}{{g(q)}\over q}\EssBdd H,
$$
\par
\noindent 
it suffices to observe that 
$$
\avesum \sum_n f(n)\wH(n-x)=\avesum \sum_{q}g(q)\sum_{{x-H\le n\le x+H}\atop {n\equiv 0(q)}}w(n-x)
=\sum_{q}g(q)\sum_{{N-H<n\le 2N+H}\atop {n\equiv 0(q)}}\sum_{{N<x\le 2N}\atop {n-H\le x\le n+H}}w(n-x)= 
$$
$$
=\sum_{q}g(q)\sum_{{n\sim N}\atop {n\equiv 0(q)}}\sum_{n-H\le x\le n+H}w(n-x) + O_{\varepsilon}\left( N^{\varepsilon}H^2\right) 
=N\FTwH(0)\sum_{q\le Q}{{g(q)}\over q} + O_{\varepsilon}\left( N^{\varepsilon}\left(H^2+QH\right)\right). 
$$
\par
\noindent 
While (ii) is a straightforward adaptation of Lemma 2.3 in [C1], in order to prove (iii) we closely follow the proof of Theorem 1.1 in [C1]. 
First, since $\Corr_{\wH}(a)$ is an even function, in (ii) we can take $K(a)=\Corr_{\wH}(a)$, $\forall a\in [-2H,2H]$, and write
$$
\sum_{a\neq 0}\Corr_{\wH}(a)R_f(a)=\sum_{\ell \le 2H}\doublesum_{(d,q)=1}g(\ell d){{g(\ell q)}\over q}\dashsum_{j(q)}\sum_{m\sim {N\over {\ell d}}}\cos {{2\pi jdm}\over q}\sum_{a}\Corr_{\wH}(a\ell)e_q(ja)+O_{\varepsilon}(N^{1+\varepsilon}H),
$$
\par
\noindent
where we have completed the last sum with the term $a=0$ at the cost
$O_{\varepsilon}(N^{1+\varepsilon}H)$ yielded by
$$
\Corr_{\wH}(0)\sum_{\ell \le 2H}\doublesum_{(d,q)=1}g(\ell d){{g(\ell q)}\over q}\dashsum_{j(q)}\sum_{m\sim {N\over {\ell d}}}\cos {{2\pi jdm}\over q}= 
$$
$$
=\Corr_{\wH}(0)\sum_{\ell \le 2H}\doublesum_{(d,q)=1}g(\ell d)g(\ell q)\Big(\sum_{{m\sim {N\over {\ell d}}}\atop {m\equiv 0(q)}}1-{1\over q}\sum_{m\sim {N\over {\ell d}}}1\Big)
\EssBdd H\sum_{\ell \le 2H}\Big(\sum_{d\le {Q\over {\ell}}}\sum_{m\sim {N\over {\ell d}}}\divisor(m)+{N\over {\ell}}\Big)  
$$
$$
\EssBdd HN\sum_{\ell \le 2H}{1\over {\ell}} 
\EssBdd HN. 
$$
\par
\noindent
Since \enspace $\widehat{\Corr}_{\wH}(\beta)\ge 0\ \forall \beta \in \R$, then from the orthogonality of the additive characters it follows that
$$
\sum_a \Corr_{\wH}(a\ell)e(a\beta)=\sum_{h\equiv 0(\ell)}\Corr_{\wH}(h)e_{\ell}(h\beta)
={1\over {\ell}}\sum_{j\le \ell}\sum_{h}\Corr_{\wH}(h)e_{\ell}(h(j+\beta))\ge 0 
\quad \forall \beta \in \R. 
$$
\par
\noindent
In particular, this yields $\displaystyle{\sum_a \Corr_{\wH}(a\ell)e_q(ja)\ge 0}\ \forall j(\bmod \; q)$, that together with the usual bound [Da, Ch.25]  implies 
$$
\sum_{a\neq 0}\Corr_{\wH}(a)R_f(a)\EssBdd \sum_{\ell \le 2H}\doublesum_{{(d,q)=1}\atop {d,q\le {Q\over {\ell}}}}{1\over q}
 \dashsum_{j(q)}{1\over {\left\Vert {{jd}\over q}\right\Vert}}\sum_{a}\Corr_{\wH}(a\ell)e_q(ja)+NH 
$$
$$
\EssBdd \sum_{\ell \le 2H}\sum_{q\le {Q\over {\ell}}}\starsum_{d\le 2q}
 {1\over q}\dashsum_{j(q)}{1\over {\left\Vert {{jd}\over q}\right\Vert}}\sum_{a}\Corr_{\wH}(a\ell)e_q(ja)+NH. 
$$
\par			
\noindent
Now, let us set $j'=jd$, so that $j=\overline{d}j'$ with \enspace $\overline{d}d\equiv 1(\bmod \, q)$, and write 
(recalling that correlations are even) 
$$
{1\over q}\dashsum_{j(q)}{1\over {\left\Vert {{jd}\over q}\right\Vert}}\sum_{a}\Corr_{\wH}(a\ell)e_q(ja)
={1\over q}\dashsum_{j'(q)}{1\over {\left\Vert {{j'}\over q}\right\Vert}}\sum_{a}\Corr_{\wH}(a\ell)e_q(j'a\overline{d})
\ll \sum_{j'\le q/2}{1\over {j'}}\sum_{a}\Corr_{\wH}(a\ell)e_q(j'a\overline{d})\ll 
$$
$$				
\ll \sum_{j\le q/2}{1\over j}\sum_{a}\Corr_{\wH}(a\ell)e_q(jan)
$$
\par
\noindent
as the variable $\overline{d}=n$ ranges over a complete set of reduced residue classes, then
$$
\starsum_{d\le 2q}{1\over q}\dashsum_{j(q)}{1\over {\left\Vert {{jd}\over q}\right\Vert}}\sum_{a}\Corr_{\wH}(a\ell)e_q(ja)
\ll \starsum_{n\le q}\sum_{j\le q/2}{1\over j}\sum_{a}\Corr_{\wH}(a\ell)e_q(jan)\ll 
$$
$$
\ll \sum_{n\le q}\sum_{j\le q/2}{1\over j}\sum_{a}\Corr_{\wH}(a\ell)e_q(jan). 
$$
\par
\noindent
Again by orthogonality of characters we get 
$$
\sum_{n\le q}\sum_{a}\Corr_{\wH}(a\ell)e_q(jan)=q\Corr_{\wH}(0)+q\sum_{{a\neq 0}\atop {ja\equiv 0(q)}}\Corr_{\wH}(a\ell), 
$$
\par
\noindent
whence 
$$
\sum_{n\le q}\sum_{j\le q/2}{1\over j}\sum_{a}\Corr_{\wH}(a\ell)e_q(jan) 
\EssBdd qH+\sum_{j\le q/2}{q\over j}\sum_{{a\neq 0}\atop {ja\equiv 0(q)}}\Corr_{\wH}(a\ell) 
\EssBdd qH+H^2, 
$$
\par
\noindent
where we have seen that
$$
\sum_{j\le q/2}{q\over j}\sum_{{a\neq 0}\atop {ja\equiv 0(q)}}\Corr_{\wH}(a\ell)
=\sum_{t'|q,t'<q}\sum_{{j\le q/2}\atop {(j,q)=t'}}{q\over j}\sum_{{a\neq 0}\atop {a\equiv 0(q/t')}}\Corr_{\wH}(a\ell)
=\sum_{t|q,t>1}\sum_{{j\le q/2}\atop {(j,q)=q/t}}{q\over j}\sum_{{a\neq 0}\atop {a\equiv 0(t)}}\Corr_{\wH}(a\ell)\ll 
$$
$$
\ll H\sum_{t|q,t\ll {H\over {\ell}}}\sum_{j'\le t/2}{t\over {j'}}\sum_{{0<|a|\le {{2H}\over {\ell}}}\atop {a\equiv 0(t)}}1 
\EssBdd H\sum_{t|q,t\ll {H\over {\ell}}}t\left({H\over {\ell t}}+1\right) 
\EssBdd H\sum_{t|q,t\ll {H\over {\ell}}}{H\over {\ell}} 
\EssBdd {{H^2}\over {\ell}}. 
$$
\par
\noindent
Finally, we have that 
$$
\sum_{a\neq 0}\Corr_{\wH}(a)R_f(a)\EssBdd NH+\sum_{\ell \le 2H}\sum_{q\le {Q\over {\ell}}}(qH+H^2)
\EssBdd NH+Q^2H\sum_{\ell \le 2H}\ell^{-2}+QH^2\sum_{\ell \le 2H}\ell^{-1}
\EssBdd NH + Q^2 H + QH^2. 
$$
\par
\noindent
The lemma is proved.\hfill $\square$ 

\bigskip

\par
\centerline{\stampatello 4. Proofs of Theorems 1, 3, 4 and Corollaries 1, 2}
\smallskip
\par
\noindent
{\bf Proof of Theorem 1}. First, let us note that
$$
\theta\in\Big({{3\vartheta-1}\over {(1-G)\vartheta+G+1}},\vartheta\Big]=
\Big(1+{{\vartheta-1}\over {\lambda_0}},\vartheta\Big] 
\Leftrightarrow
\theta = 1+{{\vartheta-1}\over {\lambda}}
$$
\par
\noindent
for some \enspace $\lambda\in(\lambda_0,1]$, where 
$$
\lambda_0\defineq 1-\vartheta+{{3\vartheta-1}\over {G+2}}>\max(1-\vartheta,1/2).
$$
\par
\noindent
Indeed, $\lambda_0>1/2$ \enspace if and only if \enspace $G+2\vartheta(1-G)>0$.
\smallskip
\par
\noindent
Then,
recalling that $\displaystyle{\chi_q(x,\sgnH)=\sum_{{{x-H}\over q}\le m\le {{x+H}\over q}}\sgn(qm-x)}$,
let us apply a dyadic argument to write
$$
J_{\sgn,g\ast \1}(N,H)=\avesum\Big|\sum_{n}\sgnH(n-x)\sum_{q|n}g(q)\Big|^2
= \avesum \Big|\sum_{q\le 2N+H}g(q)\chi_q(x,\sgnH)\Big|^2 \EssBdd 
$$
$$				
\EssBdd \max_{Q\ll N}\avesum \Big|\sum_{q\sim Q}g(q)\chi_q(x,\sgnH)\Big|^2. 
$$
\par
\noindent
From Lemma 1 (see Remark 3) one has
$$
\avesum\Big|\sum_{q\sim Q}g(q)\chi_q(x,\sgnH)\Big|^2\EssBdd NH+Q^2 H, 
\quad \forall Q\ll N. 
$$
\par
\noindent
By taking $A=Q, B=2Q$ and $Q\defineq N^{\lambda}=\infty(N/H)$ for $\lambda\in(\lambda_0,1)$ in Lemmata 2 and 3 one gets 
$$
\avesum \Big|\sum_{q\sim Q}g(q)\chi_q(x,\sgn)\Big|^2
\EssBdd \avesum\Big|\sum_{m\sim{x\over Q}}\sum_{q}g(q)\sgnH(mq-x)\Big|^2 + H^3\EssBdd 
$$
$$
\EssBdd {N\over Q}\sum_{m\asymp {N\over Q}}mJ_{\sgn,g}\Big({N\over m},\Big[{H\over m}\Big]\Big)+{{N^3}\over {Q^2}}+H^3. 
$$
\par
\noindent
Note that $QHN^{-1}=Q^\theta$ for $\theta = 1+(\vartheta-1)\lambda^{-1}$. Therefore,
$[H/m]\asymp (N/m)^{\theta}\quad \forall m\asymp N/Q$ and by hypothesis 
$$
J_{\sgn,g}\Big({N\over m},\Big[{H\over m}\Big]\Big)\ll {NH^{2-G}\over m^{3-G}} \enspace \forall m\asymp {N\over Q}. 
$$
\par
\noindent
Hence, we obtain (here, $Q=N^{\lambda}$ as above is implicit)
$$
\avesum \Big|\sum_{q\sim Q}g(q)\chi_q(x,\sgn)\Big|^2\EssBdd NH^2\Big({N\over {QH}}\Big)^G
+ {{N^3}\over {Q^2}} + H^3. 
$$
\par
\noindent
In particular, for $Q_0\defineq N^{\lambda_0}=\infty(N/H)$ it turns out that
$$
Q_0^2 H=NH^2\Big({N\over {Q_0H}}\Big)^G. 
$$
\par
\noindent
Consequently,
$$
J_{\sgn,g\ast \1}(N,H)\EssBdd NH^2\Big({N\over {Q_0H}}\Big)^G + {{N^3}\over {Q_0^2}} + H^3 
=NH^2\Big(H^{-{{G(3\vartheta-1)}\over {(G+2)\vartheta}}}+H^{-{{2(3\vartheta-1)}\over {(G+2)\vartheta}}}\Big)+H^3
\ll NH^2 \cdot H^{-{{G(3\vartheta-1)}\over {(G+2)\vartheta}}}+H^3,
$$
\par			
\noindent
that is
$$
J_{\sgn,g\ast \1}(N,H)\ll NH^{2-G'}
$$
\par
\noindent
when 
$$
0<G'<\min\left({{3-\vartheta^{-1}}\over {1+2G^{-1}}},\vartheta^{-1}-1\right). 
$$
\par
\noindent
The proof of Theorem 1 is completed.\hfill $\square$
\bigskip
\par
\noindent {\bf Proof of Theorem 3}. Recalling that $\Corr_f(0)\Corr_{\wH}(0)\EssBdd NH$,
by Lemma 4 we easily infer
$$
J_{w,f}(N,H)=
\Delta+O_{\varepsilon}(N^{\varepsilon}(NH + Q^2 H + QH^2+H^3)), 
$$
\par
\noindent
where
$$
\Delta\defineq \sum_{a\neq 0}\Corr_{\wH}(a)\sum_{\ell|a}\doublesum_{(d,q)=1}g(\ell d){{g(\ell q)}\over q}\Big( \left[ {{2N}\over {\ell d}}\right]
-\left[ {N\over {\ell d}}\right]\Big)-N\FTwH(0)^2\Big(\sum_{q\le Q}{{g(q)}\over q}\Big)^2. 
$$
\par				
\noindent
Clearly, we may confine to prove 
$$
\Delta\EssBdd NH+QH^2.
$$
\par
\noindent
To this end, observe that
$$
\sum_{a\neq 0}\Corr_{\wH}(a)\sum_{\ell|a}\doublesum_{(d,q)=1}g(\ell d){{g(\ell q)}\over q}\Big( \left[ {{2N}\over {\ell d}}\right]-\left[ {N\over {\ell d}}\right]\Big)= 
$$
$$
=N\sum_{\ell \le 2H}{1\over {\ell}}\sum_{b\neq 0}\Corr_{\wH}(\ell b)\doublesum_{(d,q)=1}{{g(\ell d)}\over d}{{g(\ell q)}\over q}
  -\sum_{a\neq 0}\Corr_{\wH}(a)\sum_{\ell|a}\doublesum_{(d,q)=1}g(\ell d){{g(\ell q)}\over q}
            \left( \left\{ {{2N}\over {\ell d}}\right\}-\left\{ {N\over {\ell d}}\right\}\right)= 
$$
$$
=N\sum_{\ell=1}^{\infty}{1\over {\ell}}\sum_{{a\neq 0}\atop {a\equiv 0(\ell)}}\Corr_{\wH}(a)\doublesum_{(d,q)=1}{{g(\ell d)}\over d}{{g(\ell q)}\over q} 
 +O_{\varepsilon}(N^{\varepsilon}QH^2), 
$$
\par
\noindent
because
$$
\sum_{a\neq 0}\Corr_{\wH}(a)\sum_{\ell|a}\doublesum_{(d,q)=1}g(\ell d){{g(\ell q)}\over q}
               \left( \left\{ {{2N}\over {\ell d}}\right\}-\left\{ {N\over {\ell d}}\right\}\right)
\ll N^{\varepsilon}H\sum_{\ell \le 2H}\sum_{0<|b|\leq {{2H}\over {\ell}}}\sum_{d\le {Q\over {\ell}}}\sum_{q\le {Q\over {\ell}}}{1\over q}
\ll N^{\varepsilon}QH^2. 
$$
\par
\noindent
Now we apply the hypothesis that $w$ is arithmetic, i.e.
$$
\sum_{{a\neq 0}\atop {a\equiv 0(\ell)}}\Corr_{\wH}(a) 
=\sum_{a\equiv 0\, (\ell)}\Corr_{\wH}(a) + O(H) 
={{\FTwH(0)^2}\over {\ell}} + O(H), 
$$
\par
\noindent
to get the conclusion 
$$
N\sum_{\ell=1}^{\infty}{1\over {\ell}}\sum_{{a\neq 0}\atop {a\equiv 0\, (\ell)}}\Corr_{\wH}(a)\doublesum_{(d,q)=1}{{g(\ell d)}\over d}{{g(\ell q)}\over q}= N\FTwH(0)^2 \sum_{\ell=1}^{\infty}{1\over {\ell^2}}\doublesum_{(d,q)=1}{{g(\ell d)}\over d}{{g(\ell q)}\over q}
  +O_{\varepsilon}\left( N^{1+\varepsilon}H\right)= 
$$
$$
= N\FTwH(0)^2 \sum_{\ell=1}^{\infty}\doublesum_{(d',q')=\ell}{{g(d')}\over {d'}}{{g(q')}\over {q'}}+O_{\varepsilon}\left( N^{1+\varepsilon}H\right)
= N\FTwH(0)^2 \Big(\sum_{n}{{g(n)}\over n}\Big)^2 + O_{\varepsilon}\left( N^{1+\varepsilon}H\right). 
$$
\par
\noindent
The proof of Theorem 3 is completed.\hfill \square
\bigskip
\par
\noindent {\bf Proof of Theorem 4}. Without further references, in what follows  we will appeal to the formula
$$
p_f\big(\log(x+m)\big)=p_f\big(\log x+\log(1+m/x)\big)
=p_f(\log x)+O\big(m(\log x)^{c-1}/x\big)
\quad \hbox{for\ $m=o(x)$ as $x\to \infty$}.
$$
\par
\noindent
After recalling that $M_f(x,H)\defineq Hp_f(\log x)$, for every $H>h$ we can write
$$
\sum_{x<n\le x+H}f(n)-M_f(x,H)
=\sum_{x<n\le x+h}f(n)-hp_f(\log x)+\sum_{x+h<n\le x+H}f(n)-(H-h)p_f\big(\log (x+h)\big)+ 
$$
$$
+H\Big(p_f\big(\log(x+h)\big)-p_f(\log x)\Big)+h\Big(p_f\big(\log x\big)-p_f\big(\log (x+h)\big)\Big). 
$$
\par
\noindent
Therefore, since $h=o(N)$ from the hypothesis $H=o(N)$, then by assuming also that $H<2h$ one has 
$$
J_f(N,H)\ll J_f(N,h) + N^{-1}H^2 h^2L^{2c-2} + \sum_{N+h<x\le 2N+h}\Big|\sum_{x<n\le x+(H-h)}f(n)-(H-h)p_f(\log  x)\Big|^2\ll 
$$
$$
\ll J_f(N,h)+N^{-1}H^2 h^2L^{2c-2} + J_f(N,H-h) + h(H-h)^2\Big(\Vert f\Vert_{\infty}^2 + L^{2c}\Big)\ll 
$$
$$				
\ll J_f(N,h) + J_f(N,H-h) + hH^2\Big(\Vert f\Vert_{\infty}^2 + L^{2c}\Big), 
$$
\par
\noindent
that gives the desired conclusion for $h<H<2h$. Now let us assume that $H\ge 2h$ and write
$$
\sum_{x<n\le x+H}f(n)=\sum_{j\le H/h}\, \sum_{x+h_{j-1}<n\le x+h_j}f(n) + \sum_{x+[H/h]h<n\le x+H}f(n), 
$$
$$
Hp_f(\log x)=\Big[ {H\over h}\Big]hp_f(\log x)+\Big\{ {H\over h}\Big\}hp_f(\log x)
=\sum_{j\le H/h}hp_f(\log x)+\Big\{ {H\over h}\Big\}hp_f(\log x) ,
$$
\par
\noindent
where we set \enspace $h_j\defineq jh$. Consequently,
$$
\sum_{x<n\le x+H}f(n)-Hp_f(\log x)=
\sum_{j\le H/h}\Big(\sum_{x+h_{j-1}<n\le x+h_j}f(n)-hp_f(\log x)\Big)+ 
$$
$$+\sum_{x+[H/h]h<n\le x+H}f(n)-\Big\{ {H\over h}\Big\}hp_f(\log x)
=\sum_{j\le H/h}\sum_{x+h_{j-1}<n\le x+h_j}f(n)-h\sum_{j\le H/h}p_f\big(\log (x+h_{j-1})\big)
+ 
$$
$$
+ \sum_{x+[H/h]h<n\le x+H}f(n)-\Big\{ {H\over h}\Big\}hp_f\big(\log (x+\big[{H\over h}\big]h)\big) 
 + O\Big({{H^2L^{c-1}}\over x}\Big). 
$$
\par
\noindent
Hence, by the Cauchy inequality 
$$
J_f(N,H)\ll {H\over h}\sum_{j\le H/h}\Big(J_f(N,h)+H\Big(\Vert f\Vert_{\infty}^2
+L^{2c}\Big)h^2\Big)+ J_f\Big(N,h\Big\{ {H\over h}\Big\}\Big) + Hh^2 \Big( \Vert f\Vert_{\infty}^2+L^{2c}\Big) + N^{-1}H^4L^{2c-2} 
$$
$$
\ll \Big({H\over h}\Big)^2 J_f(N,h) + J_f\Big(N,h\Big\{ {H\over h}\Big\}\Big) + H^3 \Big( \Vert f\Vert_{\infty}^2+L^{2c}\Big),
$$
\par
\noindent
that gives the desired conclusion also for any $H\ge 2h$.\hfill \square
\bigskip
\par
\noindent {\bf Proof of Corollary 1}. Since $d_3=\divisor\ast \1$ with $\divisor(n)\defineq\sum_{d|n}1$, we can use the bound $J_{\sgn,\divisor}(N,h)\ll Nh^{2-G}$ for every integer $h\asymp N^{\theta}$ with $\displaystyle{\theta\in(0,1/2)}$ and for every $G\in(0,1)$ (see [C-S]). Hence, from Theorem 1 we get an exponent gain $G'=G'(\vartheta,G)>0$ so that $J_{\sgn,d_3}(N,H)\ll NH^{2-G'}$ for every integer $H\asymp N^{\vartheta}$ with 
$$
\vartheta\in(1/3,1/2-\varepsilon).
$$
\smallskip
\par
\noindent
In order to exhibit a non-trivial bound for $J_{\sgn,\omega}(N,H)$, note that $\omega=\1_{\cal P}\ast \1$, where $\1_{\cal P}$ is the characteristic function of the set ${\cal P}$ of prime numbers. Thus, we can apply Theorem 1 with aid of the bound 
\smallskip
\par
\centerline{$J_{\sgn,\1_{\cal P}}(N,h)\ll Nh^{2-G}$ for every integer $h\asymp N^{\theta}$ with $\displaystyle{\theta\in\Big({1\over 6},1\Big)}$,}
\smallskip
\par
\noindent
that is a consequence of Huxley's zero density estimate [H]. Indeed, it is well known that within the same range of width $\theta\in(1/6,1)$ such an estimate implies a non-trivial bound for the classical Selberg integral (see Sect.1), whose difference from $J_{\Lambda}(N,h)$ is negligible. Thus, one obtains the aforementioned non-trivial bound for $J_{\sgn,\1_{\cal P}}(N,h)$ by using the inequalities
$$
J_{\sgn,\1_{\cal P}}(N,h)\ll J_{\1_{\cal P}}(N-h-1,h)+J_{\1_{\cal P}}(N,h),
\enspace 
J_{\1_{\cal P}}(N,h)\ll J_{\Lambda}(N,h)\log^{-2}N+Nh\log^2 N,
$$	
\par
\noindent
where the former inequality is also valid  with $\1_{\cal P}$ replaced by any $f$, while  the second one follows by applying partial summation (to $p$ powers).
\hfill \square
\bigskip
\par
\noindent {\bf Proof of Corollary 2}. Let us set $f=g\ast \1$, $\wHtwo= \CH+\uH$ and 
$$
S'_f(x)\defineq \sum_{n}\wHone(n-x)f(n),\qquad
S''_f(x)\defineq \sum_{n}\wHtwo(n-x)f(n).
$$
\par				
\noindent
Then, we apply Cauchy's inequality to write
$$
\modSel_f(N,H)-J_f(N,H)=\avesum \Big(\Big( \sum_{n}\CH(n-x)f(n)-M_f(x,H)\Big)^2-\Big( \sum_{n}\uH(n-x)f(n)-M_f(x,H)\Big)^2\Big)= 
$$
$$
=\avesum S'_f(x)\Big(S''_f(x)-2M_f(x,H)\Big)
\le\ \sqrt{J_{w',f}(N,H)\big(\modSel_f(N,H)+J_f(N,H)\big)}\ll_{\varepsilon} 
$$
$$
\ll_{\varepsilon}\sqrt{J_{w',f}(N,H)N^{1+\varepsilon}H^2}
\ll_{\varepsilon}NH^{2-3G'/2}N^{\varepsilon}\ll NH^{2-G'},
$$
\par
\noindent
after using the trivial bound for $\modSel_f(N,H)+J_f(N,H)$ and the non-trivial one yielded by Theorem 2.\hfill \square

\bigskip

\par
\centerline{\stampatello 5. Further comments and properties}
\smallskip
\par
\noindent
Because of the deep implication with the $2k-$th moments of the Riemann zeta function $\zeta$ (see [C3]), the study of the Selberg integral $J_k(N,H)$ associated to the divisor function $d_k$ has an enormous attraction. In the present section we 
add some further properties for such particular functions, in order to complement the results of [C-L].
Let us recall that in [C-L] we find the following explicit expression of the short interval mean value in the Selberg integral $J_3(N,H)$ associated to
$d_3$, i.e.
$$
M_3(x,H)\defineq H\Big( {{\log^2 x}\over 2} + 3\gamma \log x + 3\gamma^2+3\gamma_1\Big), 
$$
\par
\noindent
where $\gamma$ is the Euler-Mascheroni constant and $\gamma_1$ is a Stieltjes constant defined as 
$$
\gamma_1 \defineq \lim_{m}\Big({{\log^2 m}\over 2}-\sum_{j\le m}{{\log j}\over j}\Big).
$$
\par
\noindent
On the other side, Proposition 16 in [C-L] suggests that the expected  mean value appearing in $\modSel_3(N,H)$, where $d_3$ is weighted with $\CH$, is given by
$$
\widetilde{M}_3(x,H)\defineq H\Big(
\sum_{q\le x/N_3}{{\divisor(q)}\over q} 
 + \sum_{d_1<N_3}{1\over {d_1}}\sum_{d_2\le x/d_1 N_3}{1\over {d_2}} 
 + \Big(\sum_{d<N_3}{1\over d}\Big)^2\Big)
$$
\par
\noindent
 with $N_3\defineq[N^{1/3}]$. 
Indeed, more generally the aforementioned Proposition 16 yields the following inequality, for every $d_k$ with $k>2$ and every good weight $w$:
$$
\avesum \Big| \sum_n d_k(n)\wH(n-x)
- \FTwH(0)\sum_{q\le x/N_k}{{g_k(q)}\over q}
\Big|^2
\EssBdd 
N^{2-2/k}H+N^{-1}H^4, 
$$
\par
\noindent
where $N_k\defineq [N^{1/k}]$ and
$$
g_k(q)\defineq 
d_{k-1}(q)+\sum_{j\le k-1}\multiplesum_{{n_1 \enspace \qquad \enspace n_{k-1}}\atop {{n_1 \cdots n_{k-1} = q}\atop {n_1,\ldots,n_j<N_k}}}1
$$
\par
\noindent
is the {\it short} Eratosthenes transform of $d_k$ coming out from the $k-${\it folding} method (see [C-L], Proposition 14). 
Since Proposition 9 in [C-L] implies that
$$
M_3(x,H)-\widetilde{M}_3(x,H)
\EssBdd HN^{-1/3}
$$
\par
\noindent
uniformly for $x\sim N$,
then this justifies the presence of the {\it analytic} mean value
$M_3(x,H)$ in place of the {\it arithmetic} $\widetilde{M}_3(x,H)$ within the definition of
$\modSel_3(N,H)$. 
\par				
Here we take the opportunity to prove the next proposition, that yields for any $k>2$
the proximity in the mean square between 
the arithmetic form of the mean value 
$$
M_k(x,\wH)\defineq \FTwH(0)\sum_{q\le x/N_k}{{g_k(q)}\over q} 
$$
\par
\noindent
and its analytic counterpart
$$
\FTwH(0)p_{k-1}(\log x)=\FTwH(0)\Res_{s=1}\zeta^k(s)x^{s-1}, 
$$
\par
\noindent
where $p_{k-1}$ is the logarithmic polynomial of $d_k$. 
\medskip
\par
\noindent
To this end, first let us recall that Ivi\'c's bounds [Iv] yield an exponent gain for the Selberg integral of $d_k$,
$$
J_k(N,H)\defineq \avesum \Big|\sum_{x<n\le x+H}d_k(n)-Hp_{k-1}(\log x)\Big|^2, 
$$
\par
\noindent
once $H\asymp N^{\vartheta}$ for $\vartheta>\theta_k\defineq2\sigma_k-1$, where $\sigma_k$ is the 
{\it Carlson abscissa} for $\zeta^k$ (compare (59) of [C-L]). 
Then we prove the following property (this is Proposition
18 of [C-L], where it is stated without an explicit proof). 
\smallskip
\par
\noindent {\bf Proposition 4}.
{\it For every integer $k>2$ there exists $G=G(k)>0$ such that
$$
\avesum \left|\sum_{q\le x/N_k}{{g_k(q)}\over q}-p_{k-1}(\log x)\right|^2 \ll N^{1-G}. 
$$
} 
\smallskip
\par
\noindent {\bf Proof}. For every $H$ of width $\theta\in(0,1)$ let us take the unit step weight $\uH$ in Proposition 3. Since clearly $\FTuH(0)=H$, then Proposition 3 and the aforementioned Ivi\'c's bound imply that
$$
H^2\avesum\Big|\sum_{q\le x/N_k}{{g_k(q)}\over q}-p_{k-1}(\log x)\Big|^2 
\ll J_k(N,H) + \avesum \Big|\sum_{n}d_k(n)\uH(n-x)-\FTuH(0)\sum_{q\le x/N_k}{{g_k(q)}\over q}\Big|^2\EssBdd 
$$
$$
\EssBdd NH^2(N^{-G_1}+N^{1-2/k-\theta}+N^{-(1-\theta)}),
$$
\par
\noindent
once $\max\Big(2\sigma_k-1,1-2/k\Big)<\theta<1$. The conclusion follows for $G<\min\Big(G_1,2/k+\theta-1,1-\theta\Big)$.\hfill $\square$ 
\bigskip
\par
\noindent {\bf Remark 4}. {\it Implying in particular that $f=g\ast\1$ is a sieve function, the hypotheses of Theorem 3 legitimate the assumption that the analytic form 
$$
M_f(x,\wH)\defineq \FTwH(0)\Res_{s=1}F(s)x^{s-1}
\quad \hbox{with} \enspace F(s)\defineq \sum_{n=1}^{\infty}f(n)n^{-s}
$$
\par
\noindent
is sufficiently close to its {\it arithmetic form}, so that we can take
$$
M_f(x,\wH)= \FTwH(0)\sum_{q\le Q} {{g(q)}\over q}, 
$$
\par
\noindent
that does not depend both on $x$ and on Dirichlet series. 
From this point of view, what we have seen before shows that the divisor functions fall short of being sieve functions. Roughly speaking, they are well approximated by some arithmetic functions, which are Dirichlet convolutions of the functions $g_k$ as the $k-$folding method reveals. 
}
\medskip
\par
\noindent {\bf Acknowledgment.} This research started while the first author was a fellow \lq \lq Ing.Giorgio Schirillo\rq \rq \thinspace of the Istituto Nazionale di Alta Matematica (Italy).

\bigskip

\par				
\centerline{\stampatello References}
\smallskip
\item{[C0]} G. Coppola, On the modified Selberg integral, preprint available at http://arxiv.org/abs/1006.1229. 
\item{[C1]} G. Coppola, On the Correlations, Selberg integral and symmetry of sieve functions in short intervals,  {\it J. Combinatorics and Number Theory}, {\bf 2.2}, Article 1, 2010.
\item{[C2]} G. Coppola, On the Correlations, Selberg integral and symmetry of sieve functions in short intervals, II, {\it Int. J. Pure Appl. Math.}, {\bf 58.3}: 281--298, 2010. 
\item{[C3]} G. Coppola, On the Selberg integral of the $k$-divisor function and the $2k$-th moment of the Riemann zeta-function, {\it Publ. Inst. Math. (Beograd) (N.S.)}, {\bf 88}(102): 99--110, 2010.
\item{[C4]} G. Coppola, On the symmetry of arithmetical functions in almost all short intervals,  V, preprint available at http://arxiv.org/abs/0901.4738. 
\item{[C5]} G. Coppola, On the symmetry of square-free supported arithmetical functions in short intervals, {\it J. Inequal. Pure Appl. Math.}, {\bf 5}(2), Article 33: 1-11, 2004.
\item{[C-L]} G. Coppola and M. Laporta, Generations of correlation averages, {\it Journal of Numbers}, Vol. 2014, Article ID 140840: 1--13, 2014.
\item{[C-L1]} G. Coppola and M. Laporta, On the Correlations, Selberg integral and symmetry of sieve functions in short intervals, III,
(submitted), preprint available at http://arxiv.org/abs/1003.0302.
\item{[C-L2]} G. Coppola and M. Laporta, A generalization of Gallagher's Lemma for exponential sums, {\it \v{S}iauliai Math. Semin.}, {\bf
10} (18): 1-19, 2015.
\item{[C-L3]} G. Coppola and M. Laporta, Sieve functions in arithmetic bands, (submitted), preprint availabel at http://arxiv.org/abs/1503.07502. 
\item{[C-S]} G. Coppola and S. Salerno, On the symmetry of the divisor function in almost all short intervals, {\it Acta Arith.}, {\bf 113}(2): 189--201, 2004. 
\item{[Da]} H. Davenport, {\it Multiplicative Number Theory} Third Edition, GTM 74, Springer, New York, 2000. 
\item{[D-F-I]} W. Duke, J. Friedlander and H. Iwaniec, Bilinear forms with Kloosterman fractions, {\it Invent. Math.}, {\bf 128}(1): 23--43, 1997. 
\item{[H]} M.N. Huxley, On the difference between consecutive primes, {\it Invent. Math.}, {\bf 15}(2): 164--170, 1972.
\item{[Iv]} A. Ivi\'c, On the mean square of the divisor function in short intervals, {\it J. Th\'eor. Nombres Bordeaux}, {\bf 21}(2): 251--261, 2009. 
\item{[K-P]} J. Kaczorowski and A. Perelli, On the distribution of primes in short intervals, {\it J. Math. Soc. Japan}, {\bf 45}(3): 447--458, 1993. 
\item{[Se]} A. Selberg, On the normal density of primes in small intervals, and the difference between consecutive primes, {\it Arch. Math. Naturvid.}, 
{\bf 47}(6): 87--105, 1943. 
\smallskip
\par
\noindent
{\stampatello Giovanni Coppola}
\par
\noindent
Postal address: Via Partenio 12, 
\par
\noindent
83100 Avellino (AV), ITALY
\par
\noindent
e-mail : giocop@interfree.it
\par
\noindent
wpage: www.giovannicoppola.name
\smallskip
\par
\noindent
{\stampatello Maurizio Laporta}
\par
\noindent
Universit\`a degli Studi di Napoli "Federico II",
\par
\noindent
Dipartimento di Matematica e Applicazioni "R.Caccioppoli",
\par
\noindent
Complesso di Monte S.Angelo,
\par
\noindent
Via Cinthia, 80126 Napoli (NA), ITALY
\par
\noindent
e-mail : mlaporta@unina.it

\bye